\documentclass[12pt]{amsart}
\usepackage{a4}
\usepackage{latexsym}
\usepackage{amsmath}
\usepackage{amssymb}
\usepackage{amscd}
\usepackage[all]{xy}
\usepackage{paralist}
\usepackage{color}
\usepackage{amsthm}

\theoremstyle{plain}
    \newtheorem{theorem}                    {Theorem}[section]  
    \newtheorem{lemma}      [theorem]       {Lemma}
    \newtheorem{corollary}  [theorem]       {Corollary}
    \newtheorem{proposition}[theorem]       {Proposition}

\newtheorem{example}[theorem]{Example}

\newcommand{\X}{{\mathcal X}}

\newcommand{\XK}{X}

\DeclareFontFamily{U}{wncy}{}
    \DeclareFontShape{U}{wncy}{m}{n}{<->wncyr10}{}
    \DeclareSymbolFont{mcy}{U}{wncy}{m}{n}
    \DeclareMathSymbol{\Sh}{\mathord}{mcy}{"58}

\newcommand{\Hom}{\operatorname{Hom}}
\newcommand{\End}{\operatorname{End}}
\newcommand{\Ext}{\operatorname{Ext}}

\newcommand{\Pic}{\operatorname{Pic}}

\newcommand{\tor}{{\operatorname{tor}}}
\newcommand{\Tor}{{\operatorname{Tor}}}

\newcommand{\cont}{{\operatorname{cont}}}

\newcommand{\Spec}{\operatorname{Spec}}

\newcommand{\im}{\operatorname{im}}

\newcommand{\coker}{\operatorname{coker}}
\newcommand{\id}{\operatorname{id}}
\newcommand{\Br}{\operatorname{Br}}
\newcommand{\CH}{\operatorname{CH}}
\renewcommand{\lim}{\operatorname{lim}}
\newcommand{\colim}{\operatornamewithlimits{colim}}
\newcommand{\rk}{\operatorname{rank}}

\newcommand{\NS}{\operatorname{NS}}

\newcommand{\Z}{{{\mathbb Z}}}
\newcommand{\Zp}{{\mathbb Z}_{(p)}}
\newcommand{\Q}{{{\mathbb Q}}}

\newcommand{\F}{{{\mathbb F}}}
\renewcommand{\L}{{{\mathbb L}}}
\newcommand{\G}{{{\mathbb G}}}

\newcommand{\K}{{\mathcal K}}

\newcommand{\et}{{\text{\rm et}}}

\newcommand{\cone}{\operatorname{cone}}
\renewcommand{\O}{{\mathcal O}}

\newcommand{\proofend}{\hfill$\square$\\ \smallskip}

\date{\today}
\title{On the kernel of the Brauer-Manin pairing}
\author{Thomas H. Geisser}
\address{
	Department of Mathematics, Rikkyo University, Ikebukuro, Tokyo, Japan
}
\email{geisser@rikkyo.ac.jp}
\author{Baptiste Morin}
\address{
	Department of Mathematics, Universit\'e de Bordeaux,
	Bordeaux, France
}
\email{Baptiste.Morin@math.u-bordeaux.fr}
\thanks{The first named author is Supported by JSPS Grant-in-Aid (C) 18K03258,
and the second named author by grant ANR-15-CE40-0002}
\subjclass[2010]{Primary:\ 14F22;\ Secondary:\ 14F40, 14F42, 11G40}
\keywords{Brauer group; Brauer Manin pairing; local fields; Picard number}

\begin{document}

\begin{abstract}
Let $\X$ be a regular scheme, flat and proper over the ring of integers
of a $p$-adic field, with generic fiber $X$ and special fiber $\X_s$.
We study the left kernel $\Br(\X)$ of the Brauer-Manin pairing
$\Br(X)\times \CH_0(X)\to \Q/\Z$. Our main result is that
the kernel of the reduction map $\Br(\X)\to \Br(\X_s)$ is the direct sum
of $(\Q/\Z[\frac{1}{p}])^s\oplus (\Q/\Z)^t$ and a finite $p$-group, 
where $s+t= \rho_{\X_s}-\rho_X-I+1$, for $\rho_{\X_s}$ and $\rho_X$
the Picard numbers of $\X_s$ and $X$, and $I$ the number of irreducible
components of $\X_s$. Moreover, we show that $t>0$ implies $s>0$.
\end{abstract}

\maketitle

\section{Introduction}
The Brauer group plays an important role in arithmetic geometry.
Over a finite field, Artin conjectured that the Brauer group
of any proper scheme is finite \cite[Rem. 2.5c)]{brauerIII}; this was
proved by Grothendieck for curves %\cite[Cor. 1.2, Rem. 2.5b)]{brauerIII}.
\cite[Rem. 2.5b)]{brauerIII}. If $X$ is smooth and proper,
then the finiteness of $\Br(X)$ is equivalent to Tate's conjecture 
on the surjectivity of the cycle map 
for divisors on $X$, and for a normal crossing scheme the finiteness 
follows from Tate's conjecture for all (smooth) intersections of the 
components. The next interesting case are varieties over a $p$-adic field $K$.
It is a classical result of Hasse that the Brauer group $\Br(K)$
is isomorphic to $\Q/\Z$. 
%and the Brauer group of
%a global field is contained in a short exact sequence
%$$ 0\to \Br(K)\to\oplus_v\Br(K_v)\stackrel{\Sigma}{\to}\Q/\Z\to 0.$$
Lichtenbaum  \cite{licht} proved that if $X$ is a curve, then the Brauer group 
$\Br(X)$ is Pontrjagin dual to the Chow group of zero cycles
$\CH_0(X)$. In particular, it is the direct sum
of a finite group, of $\Q/\Z$, and of a divisible $p$-torsion group
of corank the  genus of $X$ times the degree of $[K:\Q_p$].

This result was generalized by Colliot-Th\'el\`ene and Saito \cite{ctsaito}, 
and Saito and Sato \cite{saisa}. 
They show that if $X$ has a proper and regular model $\X$, then the Brauer-Manin
pairing between $\CH_0(X)$ and $\Br(X)$ has left kernel $\Br(\X)$.
Moreover, for $l\not=p$, the $l$-part of $\Br(X)/\im(\Br(\X)\oplus\Br(K))$
is finite and vanishes for almost all $l$. 
However, not much is known about $\Br(\X)$. We prove the following:

\begin{theorem}\label{maintheo}
Let $\X$ be a regular scheme, flat and proper over the ring of integers
of a $p$-adic field, and let $\X_s$ be the closed fiber. Then 
the kernel of the reduction map $\Br(\X)\to \Br(\X_s)$ is the direct sum
of $(\Q/\Z[\frac{1}{p}])^s\oplus (\Q/\Z)^t$ and a finite $p$-group, 
where $s+t=r:=\rho_{\X_s}-\rho_X-I+1$ for $\rho_{\X_s}$ and $\rho_X$
the Picard numbers of $\X_s$ and $X$, and $I$ the number of irreducible
components of $\X_s$. Moreover, if $t>0$ then $s>0$.
\end{theorem}

Note that $\Br(\X_s)$ is conjecturally finite.
The statement on the $l$-corank follows from the proper
base change theorem, and the non-trivial part of the theorem is that 
the $p$-corank is strictly smaller than 
the $l$-corank unless both vanish.

\begin{corollary}
1) The kernel of $\Br(\X)\to \Br(\X_s)$ is finite if and only if 
$r=0$. 

2) If $r=1$, then the $p$-part of the kernel of 
$\Br(\X)\to \Br(\X_s)$ is finite.
\end{corollary}

%Using a theorem of Flach-Siebel, we relate the rank of the
%$p$-part of the Brauer group to $h^{0,2}=\dim H^2(X,\O_X)$.
Our construction together with a theorem of Flach-Siebel
gives a map $b: \Pic\X_s\to H^2(X,\O_X)$ which
is related to the Chern class map. We show in Theorem \ref{chern} 
that $s$ is the dimension of the $\Q_p$-vector 
space spanned by image of $\Pic(\X_s)$ in $H^2(X,\O_X)$.
In particular, $H^2(X,\O_X)=0$ implies that the kernel of 
$\Br(\X)\to \Br(\X_s)$ is finite. We use this to give some
explicit calculations.

\begin{theorem}
Let $\X$ be a family of abelian or K3-surfaces over $\Spec \Z_p$.
If $r=0$, then $\Br(\X)$ is finite. If $r>0$, then 
$$\textstyle
\Br(\X)\cong (\Q/\Z[\frac{1}{p}])\oplus (\Q/\Z)^{r-1}
\oplus P,$$
where $P$ is a finite $p$-group.
\end{theorem}

We give an explicit example of an abelian surface with
$$\textstyle
\Br(\X)\cong (\Q/\Z[\frac{1}{p}])\oplus (\Q/\Z)^2\oplus P.$$
Finally, we briefly discus the intermediate groups
$$ \Br(\X)\to \lim \Br(\X_n)\to \Br(\X_s)$$
for $\X_n=\X\times_\Z\Z/p^n\Z$.

\subsection*{Notation:}
Throughout the paper, 
$K$ is a finite extension of $\Q_p$ of degree $f$ with Galois group $G_K$, 
and $X$ be a smooth and proper scheme over $K$ of dimension $d$.
We let $h^{0,i}=\dim_K H^i(X,\O_X)$, and $\rho_X=\rk \NS(X)$
the Picard number. 

We let $\O_K$ be the ring of integers of $K$ and
assume that $X$ has a proper regular model $\X/\O_K$, which we can 
(by Stein factorization) assume to have geometrically connected fibers.
Let $i:\X_s\to \X$ be the special fiber, $\rho_{\X_s}=\rk \Pic(\X_s)$ 
its Picard number and $I$ the number of irreducible components of $\X_s$.
%We define the Picard number of $\X$ to be $\rho_\X=\rho_X+I-1$.
% (see Proposition\ref{saisa}(2)). 
The number 
$$r= \rho_{\X_s}-\rho_X-I+1$$ 
plays an important role in this paper. 

For an abelian group $A$ we let $A^{\wedge l}=\lim_i A/l^i$ be
the $l$-adic completion, ${}_mA$ be the subgroup of $m$-torsion 
elements, $T_l=\lim_i {}_{l^i}A$ the $l$-adic Tate module,
and $V_l=T_l\otimes_{\Z_l}\Q_l$.

\subsection*{Acknowledgements:}
We would like to thank T. Suzuki for helpful discussions. The second
named author thanks Rikkyo University for their hospitality.

\section{The Brauer group}
We start by recalling some known facts on the cohomology of $\G_m$. 
Recall that $f=[K:\Q_p]$.

\begin{proposition}\label{saisa}
1) We have $H^0_\et(\X,\G_m)\cong \O(\X)^\times$, a direct sum of a 
$\Z_p$-module of rank $f$ and a finite group, and 
$H^0_\et(\XK,\G_m)\cong \O(\XK)^\times\cong \O(\X)^\times\times \Z$.

2) The group $H^1_\et(\XK,\G_m)\cong \Pic(\XK)$ is an 
extension of a finitely generated group of rank $\rho_X$ 
by a finitely generated $\Z_p$-module of rank $f\cdot h^{0,1}$, 
and there is an exact sequence
\begin{equation}\label{picseq}
% 0 \to\O^\times\to  K^\times \to \Z^I\\
0\to \Z^{I-1}\to \Pic(\X) \to \Pic(X)\to 0,
\end{equation}
%$$ \oplus_{\X_s^{(0)}}\Z\to \Pic(\X)\to \Pic(\XK)\to 0.$$
where $I$ is the number of irreducible components of $\X_s$.

3) The Brauer groups $H^2_\et(\X,\G_m)\cong \Br(\X)$ and 
$H^2_\et(\XK,\G_m)\cong \Br(\XK)$ are torsion groups with finite 
$m$-torsion for every $m$.
Moreover, $\Br(\XK)$ contains $\Br(\X)$ and $\Q/\Z\cong\im \Br(K)$
as subgroups with trivial intersection, 
and $\Br(\XK)/(\im \Br(K)\oplus \Br(\X))$ is isomorphic
to the sum of a finite group and finitely many copies of $\Q_p/\Z_p$.
\end{proposition}

\proof
1) This follows from $\O(\X)\cong \O_K$ and $\O(\XK)\cong K$
because of geometric connectedness.

2) %Since the induced map on $\Pic$ is injective with cokernel of finite
%index for finite field
%extensions we can assume that $\XK$ has a $K$-rational point.
Consider the low degrees term of the Hochschild-Serre
spectral sequence:
$$ 0\to \Pic(X)\to \Pic(\bar X)^{G_K} \stackrel{d_2}{\to}\Br(K)\to \Br(X).$$
Since $\Q/\Z\cong \Br(K)\to \Br(X)$ has finite kernel (as one sees 
with a $K'$-rational point for $K'/K$ finite),
the image of $\Br(K)$ in $\Br(X)$ is isomorphic to $\Q/\Z$
and $\Pic(X)$ and  $\Pic(\bar X)^{G_K}$ differ by a finite group.
 
We know that $\Pic(\bar \XK)^{G_K}$ is %up to a finite group
an extension of the finitely generated
N\'eron-Severi group (of rank $\rho_X$) and the rational points of 
an abelian variety of dimension $g$, which has a subgroup of finite
index isomorphic to $\O_K^g$ by Mattucks's theorem \cite{mattuck}. 
In view of $H^i_\et(\X,\G_m)\cong CH^1(\X,1-i)$ and 
$H^i_\et(X,\G_m)\cong CH^1(X,1-i)$ for $i\leq 1$,  
the sequence is the localization sequence for higher Chow groups
$$ 0 \to\O^\times\to  K^\times \to \Z^I
\to \Pic(\X) \to \Pic(X)\to 0,$$
where we use the identification $CH^1(\X,1)=\O^\times$, 
$CH^1(X,1)=K^\times$, as well as $CH_d(\X_s)\cong \Z^I$,
the free abelian group on the irreducible components of $\X_s$.

3) The groups are torsion groups because Brauer groups of regular schemes
are contained in the Brauer group of their function field. 
Since $\Br(\X)\subseteq \Br(X)$, it suffices to prove
finiteness of the $m$-torsion for $\Br(\XK)$. 
It suffices to prove that $H^2_\et(X,\mu_m)$ is finite, because it 
surjects onto ${}_m\Br(X)$. This finiteness follows from
the Hochschild-Serre spectral sequence
$$H^s(K,H^t_\et(\bar X,\mu_m))\Rightarrow H^{s+t}_\et(X,\mu_m)$$
because the coefficients $H^t_\et(\bar X,\mu_m)$ are finite, and 
Galois cohomology of a local field of characteristic $0$ of a finite 
module is finite.

The prime to $p$-part of the statement about 
$\Br(\XK)/(\im \Br(K)\oplus \Br(\X))$ is proven in 
\cite{ctsaito}, see also \cite[Prop. 5.2.1]{saisa}. The $p$-part
follows from the finiteness of the $p$-torsion of $\Br(\XK)$.
\proofend

\medskip 

For later use we note the following facts about 
the cohomology of the special fiber. 

\begin{proposition}
The group of units $H^0_\et(\X_s,\G_m)$ is a finite group, the Picard group
$H^1_\et(\X_s,\G_m)$ is finitely generated, and the Brauer 
group is the torsion subgroup of $H^2_\et(\X_s,\G_m)$.
\end{proposition}

The first statement is clear, the second statement 
can be found in \cite{roquette}, and the third statement is
an unpublished theorem of Gabber \cite{g-dejong}.

\subsection*{Estimates using $l$-adic cohomology}
For any prime $l$ (including $l=p$), we have the short exact
coefficient sequence 
$$ 0\to \Pic(X)^{\wedge l}\otimes_{\Z_l}\Q_l \to H^2(X,\Q_l(1))\to 
V_l\Br(X) \to 0.$$
The left $\Q_l$-vector space has dimension equal rank $\rho_X$
if $l\not=p$, and equal to $\rho_X+fh^{0,1}$ for $l=p$ by Proposition 
\ref{saisa}(2). %because $h^{0,1}$ is the dimension of $\Pic^0_X$.
Thus in order to understand $V_l\Br(X)$, we calculate $H^2(X,\Q_l(1))$.
The spectral sequence
\begin{equation}\label{hochserre}
E_2^{s,t}= H^s(K,H^t(\bar X,\Q_l(1)))\Rightarrow H^{s+t}(X,\Q_l(1))
\end{equation}
degenerates at $E_2^{s,t}$ by \cite{deligne}. Since  
$H^2(K,H^0(\bar X,\Q_l(1)))\cong V_l\Br(K)\cong \Q_l$, we obtain 
$$\dim H^2(X,\Q_l(1))= 1+ \dim  H^1(K,H^1(\bar X,\Q_l(1)))
+\dim H^2(\bar X,\Q_l(1))^{G_K}.$$
From the divisibility of $\Pic^0(\bar X)$ we obtain a short exact sequence
$$ 0\to \NS(\bar X)\otimes\Q_l \to 
H^2(\bar X,\Q_l(1))\to V_l\Br(\bar X)\to 0.$$
The vanishing of $H^1(K,\NS(\bar X)\otimes\Q_l)$ then implies
$$\dim H^2(\bar X,\Q_l(1))^{G_K}= 
\dim \NS(X)_{\Q_l}+\dim  V_l\Br(\bar X)^{G_K}=
\rho_X+\dim V_l\Br(\bar X)^{G_K}.$$
%Next we consider the vector space $H^1(K,H^1(\bar X,\Q_l(1)))$. 
The remaining direct summand of $H^2(X,\Q_l(1))$
is calculated in the following proposition.

\begin{proposition}
The vector space $H^1(K,H^1(\bar X,\Q_l(1)))$ vanishes if $l\not=p$,
and it has dimension $2fh^{0,1}$ if $l=p$.
\end{proposition}

\proof 
We have $H^1(\bar X,\Q_l(1))\cong V_l\Pic^0(\bar X)$ is a vector space of dimension 
$2h^{0,1}$ for any $l$. Using Euler-Poincar\'e characteristic,
$$
\prod_i H^i(G_K,V)^{(-1)^i}=
\begin{cases} 
1, & V\ \text{a $\Q_l$-vector space};\\
-f\dim V, & V\ \text{a $\Q_p$-vector space},\\
\end{cases}$$
it suffices to show that 
$H^1(\bar X,\Q_l(1))^{G_K}$ and $H^2(K, H^1(\bar X,\Q_l(1)))$ vanish.
To show the vanishing of $H^1(\bar X,\Q_l(1))^{G_K}$, 
we note that $V_l\Pic(X)=0$ implies that the two left groups in the 
short exact sequence arising from \eqref{hochserre},
$$ 0\to H^1(K,\Q_l(1))\to H^1(X,\Q_l(1))\to 
H^1(\bar X,\Q_l(1))^{G_K}\to 0$$
are both isomorphic to $(K^\times)^{\wedge l}\otimes_{\Z_l}\Q_l$. 
By local duality  
$$H^2(K, H^1(\bar X,\Q_l(1))\cong H^1(\bar X,\Q_l)_{G_K}$$
and by Poincar\'e duality the right hand term is dual to 
$H^{2d-1}(\bar X,\Q_l(d))^{G_K}$, which, by the hard Lefschetz 
theorem, is isomorphic to $H^{1}(\bar X,\Q_l(1))^{G_K}=0$.
\qed

Comparing the two expressions for $H^2(X,\Q_l(1))$, we obtain

\begin{theorem}
We have
$$\dim_{\Q_l} V_l\Br(X)=\begin{cases}
1+\dim V_l\Br(\bar X)^{G_K}&\text{for } l\not=p;\\
1+\dim V_p\Br(\bar X)^{G_K}+fh^{0,1} &\text{for } l=p.
\end{cases}$$
\end{theorem}

Note that in general 
$\dim_{\Q_p} H^2(\bar X,\Q_p(1))^{G_K}\leq\dim_{\Q_l} 
H^2(\bar X,\Q_l(1))^{G_K}$,
or equivalently,
$\dim_{\Q_p} V_p\Br(\bar X)^{G_K} \leq\dim_{\Q_l}V_l\Br(\bar X)^{G_K}$
for $l\not=p$,  so that it is not clear which of 
$\dim_{\Q_l} V_l\Br(X)$ and $\dim_{\Q_p} V_p\Br(X)$ is larger.

If $\X$ is a regular proper model, then
%$\rho_X=\rho_\X-I+1$, for $I$ the number of irreducible 
%components of the special fiber $\X_s$.
by Proposition \ref{saisa}(3) and the proper base change theorem,
we have  
$$\dim_{\Q_l} V_l\Br(X)=1+ \dim_{\Q_l} V_l\Br(\X)=1+r
+\dim_{\Q_l} V_l\Br(\X_s)$$ for $l\not=p$, 
hence 
%\begin{align*}
$$\dim V_l\Br(\bar X)^{G_K}=r+\dim_{\Q_l} V_l \Br(\X_s)%\\
%\dim H^2(\bar X,\Q_l(1))^{G_K}& =\rho_s-I+1+\dim_{\Q_l}V_l \Br(\X_s).
$$%\end{align*}
In particular, $r$ is independent of
the model if we assume Artin's conjecture on the finiteness of $\Br(\X_s)$.

\section{Completions}\label{completes}
We five some facts about completions needed below;
the reader familiar with completions can skip this section.
For a complex $A$ of abelian groups and integer $m$, 
$A\otimes^\L\Z/m$ is represented by
the total complex of the double complex $A\stackrel{m}{\to} A$
concentrated in (cohomological) degrees $-1$ and $0$. 
The canonical map $A\to A\otimes^\L\Z/m$ is induced by mapping $A$
to the component in degree $0$.
The cohomology of $A\otimes^\L\Z/m$ can be calculated by the exact sequence 
\begin{equation}\label{coeffs}
 0\to H^i(A)/m \to H^i(A\otimes^\L\Z/m)\to {}_m H^{i+1}(A) \to 0.
\end{equation}
For a prime number $p$, the $p$-completion of $A^\wedge$ is the 
pro-system $\{A\otimes^\L\Z/p^j\}_j$,
where the transition maps in the system are multiplication
by $p$ in degree $-1$ and the identity in degree $0$. 
We define continuous cohomology $H^i_\cont(A^\wedge)$ to be the cohomology of
$R\lim A^\wedge$. We have a short exact sequence
$$ 0\to \lim^1_j H^{i-1}(A\otimes^\L\Z/p^j) \to H^i_\cont(A^\wedge)
\to \lim_j H^{i}(A\otimes^\L\Z/p^j)\to 0.$$
With the inverse limit of the sequence \eqref{coeffs} and using 
$\lim^1_j H^i(A\otimes^\L\Z/p^j)\cong \lim^1_j {}_{p^j}H^{i+1}(A)$ we obtain a 
diagram with exact rows and columns
$$\begin{CD}
@.@.@.0\\
@.@.@.@VVV\\
@.@. H^i(A)@>d^i>> \lim_j H^i(A)/p^j\\
@.@.@Vc^iVV @VVV \\
0@>>>\lim^1_j {}_{p^j}H^i(A)@>>> H^i_\cont(A^\wedge)
@>e^i>> \lim_j H^i(A\otimes^\L\Z/p^j)@>>> 0\\
@.@.@.@VVV\\
@.@.@. T_p H^{i+1}(A)\\
@.@.@.@VVV\\
@.@.@.0
\end{CD}.$$
If ${}_{p^j}H^i(A)$ is finite for one (or, equivalently, for all) 
$j$, then $\lim^1_j {}_{p^j}H^i(A)$
%and does not contain any divisible subgroup \cite[Lemma 4.3c)]{jannsen}.
vanishes and $\ker c^i=\ker d^i$.
%inverse limit of the sequence \eqref{coeffs} gives the sequence
%$$ 0\to H^i(A)^{\wedge p} \to H^i_\cont(A^\wedge)\to 
%T_p H^{i+1}(A)\to 0.$$
The kernel of $d^i$ consists of the $p$-divisible elements of $H^i(A)$ but we
are interested in the (smaller) kernel of $c^i$.

\begin{proposition}\label{compllemma}
The kernel of $c^i: H^i(A)\to H^i_\cont(A^\wedge)$ is the maximal $p$-divisible
subgroup of $H^i(A)$. The cokernel of $c^i$ is an extension of  
$T_pH^{i+1}(A)$ by the uniquely $p$-divisible group 
$\lim^1 (H^i(A),p)$. In particular, it is $p$-torsion free.
%, and its torsion subgroup is exactly 
%the maximal $p$-divisible torsion subgroup of $H^i(A)$.
\end{proposition}

%Note that the kernel may be strictly larger than the maximal divisible subgroup.
%For example, the group 
%$\oplus_{i=0}^{\infty} \Z e_i/\langle p^ie_i=e_0\rangle$ is 
%reduced, but the kernel is isomorphic to $\Z$, generated by the $p$-divisible
%element $e_0$.

We will give some calculations of $\lim^1 (H^i(A),p)$
below in case $H^i(A)$ are finitely generated of rank $r$.
% then 
%$\lim^1 (H^i(A),p)\cong\Q^\aleph\oplus (\Q/\Z')^r $
%with $\aleph$ the cardinality of the continuum.

\proof
The diagram gives an exact sequence
\begin{equation}\label{666}
0\to \ker c^i \to \ker d^i \to \lim^1_j {}_{p^j}H^i(A)
\to \coker c^i\to \coker e^ic^i\to 0.
\end{equation}
For any abelian group $G$ we have a sequence of inverse systems
$$\begin{CD}
0@>>> {}_{p^{j+1}}G @>>> G @>p^{j+1}>> p^{j+1}G@>>> 0\\
@.@VVV @VpVV @VVV \\
0@>>> {}_{p^j}G @>>> G @>p^j>> p^jG@>>> 0\\
\end{CD}
$$
and taking the $6$-term exact derived lim sequence we get
\begin{equation}\label{66}
0\to T_pG\to \lim (G,p)\stackrel{\xi}{\to} \lim p^jG\to \lim^1 {}_{p^j}G\to 
\lim^1 (G,p)\to \lim^1 p^jG\to 0.
\end{equation}
Applying this to $G=H^i(A)$ and comparing to \eqref{666} we see that 
$\lim p^j H^i(A)= \bigcap_j p^j H^i(A)=\ker d^i$ implies 
$\ker c^i\cong \im\xi$, which is the maximal 
$p$-divisible group of $G$ by \cite[Lemma 4.3a)]{jannsen}. 

By definition, the cone of the completion map $c$ of complexes is isomorphic
to the cohomology of $R\lim (A,p)$, where the system is the complex
$A$ with transition map multiplication by $p$. We obtain a diagram
$$\begin{CD}
0@>>> \lim^1 (H^i(A),p)@>>> H^i(\cone(c)) @>>> \lim (H^{i+1}(A),p)@>>> 0\\
@.@VVV @|@V\xi VV \\
0@>>> \coker c^i@>>> H^i(\cone(c)) @>>> \ker c^{i+1}@>>> 0\\
\end{CD}$$
By \eqref{66}, the kernel of the right map is 
$T_pH^{i+1}(A)$ and we see by the snake Lemma that $\coker c^i$ is an
extension of $T_pH^{i+1}(A)$ by $\lim^1 (H^i(A),p)$.
\proofend

\begin{corollary}
If $H^i(A)$ is a $\Zp$-module, 
then $\coker c^i$  is the direct sum of $T_pH^{i+1}(A)$ and the 
uniquely divisible group $\lim^1 (H^i(A),p)$.
\end{corollary}

A similar discussion applies to a bounded complex of sheaves $B$ on $X$ by 
applying the above to $A=R\Gamma_\et(X,B)$. Since
$R\Gamma_\et(X,B\otimes^\L\Z/m)\cong R\Gamma_\et(X,B)\otimes^\L\Z/m$,
%and $R\lim R\Gamma_\et(X,-)=R\Gamma_\et(X,R\lim -)$
we define $H^i_\cont(X,B^\wedge)$ as the cohomology of 
$R\lim_j (R\Gamma_\et(X,B)\otimes^\L\Z/p^j)$
%\cong R\Gamma_\et(X,R\lim(B\otimes^\L\Z/p^r)).\end{equation}
and obtain the sequence
\begin{equation}\label{contcohom}
0\to H^i_\et(X,B)^{\wedge p} \to H^i_\cont(X,B^\wedge)\to 
T_p H^{i+1}_\et(X,B)\to 0
\end{equation}
if $H^{i-1}_\et(X,B\otimes^\L\Z/p)$ is finite.

\section{The $p$-adic completion of $\X$}
We have an exact sequence of etale sheaves on $\X$, 
\begin{equation}\label{sheafes}
0\to \K \to \G_{m,\X} \to i_* \G_{m,\X_s}\to 0.
\end{equation}
Since $\mu_m=\ker (\G_m\stackrel{\times m}{\to} \G_m)$ is locally constant 
on $\X$ for $m$ prime to $p$, the 
proper base change theorem implies that the cohomology of $\K$ is uniquely
$l$-divisible for any $l\not=p$, i.e., it consists of $\Zp$-modules.

As the map  $\O(\X)^\times\to \O(\X_s)^\times$ is surjective
with $p$-adically complete kernel $H^0(\X,\K)$, we
obtain a long exact sequence
\begin{equation}\label{seq}
0\to H^1_\et(\X,\K) \to \Pic(\X)\to\Pic(\X_s) \to H^2_\et(\X,\K)
\to \Br(\X)\to \Br(\X_s).
\end{equation}
The group $\Pic(\X_s)$ is finitely generated \cite{roquette}, and we
let $N\subseteq H^2_\et(\X,\K)$ be the cokernel of the map 
$\Pic(\X)\to \Pic(\X_s)$.

\begin{proposition}\label{KandN}
The group $N$ is finitely generated of rank $r= \rho_{\X_s}-\rho_X-I+1$,
and $H^1_\et(\X,\K)$ is a finitely generated $\Z_p$-module of rank 
$f h^{0,1}$.
\end{proposition}

\proof
Since $N$ is finitely generated it suffices to calculate
$N\otimes \Q_l/\Z_l$ for $l\not=p$. We have a short exact sequence 
$$0\to \Pic(\X)\otimes \Q_l/\Z_l\to \Pic(\X_s)\otimes 
\Q_l/\Z_l\to N\otimes \Q_l/\Z_l\to 0,$$
where the left map is injective by the proper base change theorem. 
The Lemma follows from Proposition \ref{saisa}(2) by counting coranks.

The same short exact sequence shows that
the subgroup $\Z^{I-1}\subseteq \Pic(\X)$ injects into $\Pic(\X_s)$.
Hence $H^1_\et(\X,\K)$ can be viewed a subgroup of $\Pic(X)$. Since it
is uniquely $l$-divisible, it maps to zero in $\NS(X)/l$ for all
$l\not=p$, hence it has trivial image in the finitely generated
group $\NS(X)$. Thus 
$H^1_\et(\X,\K)$ can be viewed as a subgroup of $\Pic^0(X)$. 
Since the quotient is finitely generated, we conclude that it is a 
finitely generated $\Z_p$-module of the same rank $fh^{0,1}$. 
\endproof

The $\Zp$-module $H^2_\et(\X,\K)$ has finite $p$-torsion because $N$ is 
finitely generated and $\Br(\X)$ has finite $p$-torsion.
Since $\Br(\X)$ is torsion, $N_\Q\cong H^2_\et(\X,\K)_\Q$ has dimension $r$.

\medskip

Consider the formal completion of $\X$  at $p$, i.e., the
direct system of the reductions $\X_n=\X\times_\O\O/p^n$, 
of $\X$ modulo $p^n$. 
We obtain a short exact sequence of
pro-sheaves on the topological space $\X_s$,
\begin{equation}\label{complcoh}
0\to \K_n \to i_n^*\G_{m,\X_n} \to \G_{m,\X_s}\to 0,
\end{equation}
where $i_n:\X_s\to \X_n$ is the closed embedding. 
Let $H^i_\cont(\X_s, \G_{m,\bullet})$  and $H^i_\cont(\X_s, \K_\bullet)$
be the continuous cohomology of the 
pro-sheaves $(i_n^*\G_{m,\X_n})_n$ and $(\K_n)_n$ on $\X_s$, respectively.

Since $i_n$ is a universal homeomorphism, 
the cohomology can be calculated with the short exact sequence
\begin{equation}\label{slim}
0\to \lim^1 H^{i-1}_\et(\X_n,\G_{m,\X_n}) \to H^i_\cont(\X_s, \G_{m,\bullet})
\to \lim H^{i}_\et(\X_n,\G_{m,\X_n})\to 0.
\end{equation}
If $i=1$, then the left term vanishes because the groups 
$H^0_\et(\X_n,\G_{m,\X_n})$ are finite by properness of $\X_n$. Moreover, 
the natural map $\Pic(\X)\to \lim\Pic(\X_n)$ is an isomorphism by 
Grothendiecks formal existence theorem 
\cite[Cor. 5.1.6, Scholie 5.1.7]{EGAIII}, so that we have an isomorphism 
\begin{equation}\label{piciso}
\Pic(\X)\cong  H^1_\cont(\X_s, \G_{m,\bullet}).
\end{equation}
This implies %Using \eqref{piciso}, 
that the long exact cohomology sequence associated
to the short exact sequence \eqref{complcoh} takes the form
\begin{multline}\label{seqc}
0\to H^1_\cont(\X_s,\K_\bullet) \to \Pic(\X) \to\Pic(\X_s)\\ 
\to H^2_\cont(\X_s,\K_\bullet)\to H^2_\cont(\X_s,\G_{m,\bullet}) \to 
H^2_\et(\X_s,\G_m).
\end{multline}

Let $H^i_\cont(\X,\K^\wedge)$ be the cohomology of the $p$-adic completion 
$\K^\wedge$ of $\K$ as in Section \ref{completes}.

\begin{proposition}\label{mainth}
We have $H^i_\cont(\X_s,\K_\bullet) \cong H^i_\cont(\X,\K^\wedge)$ 
for all $i$. 
%$$\lim H^2(\X,\K)/p^r \cong \hat H^2(\X,\K)\cong H^2_\cont(\X,\Kc).$$
\end{proposition}

To prove the proposition, we compare both sides to the $p$-completion 
$(\K_n^\wedge)_n$, i.e., the double inverse 
system $i_n^*\K_n\otimes^\L\Z/p^t$, and show that we have isomorphisms
$$ H^i_\cont(\X,\K^\wedge)\stackrel{\sim}{\longrightarrow}  
H^i_\cont(\X_s,i^*\K^\wedge)\stackrel{\sim}{\longrightarrow}
H^i_\cont(\X_s,\K^\wedge_\bullet)\stackrel{\sim}{\longleftarrow} 
H^i_\cont(\X_s,\K_\bullet). $$
The first isomorphism follows from the proper base change because 
$H^i(\O_K,\mathcal F)\stackrel{\sim}{\longrightarrow} H^i(s,i^*\mathcal F)$ 
for any etale sheaf $\mathcal F$ on $\Spec \O_K$.
The third and second isomorphism follow from the following Proposition.

\begin{proposition}
1) For fixed $t$, we have an isomorphism of pro-sheaves on $\X_s$
$$ i^*\K\otimes^\L\Z/p^t \stackrel{\sim}{\to} \{i_n^*\K_n\otimes^\L\Z/p^t\}_n.$$

2) For fixed $n$, we have an isomorphism of pro-sheaves on $\X_n$,
$$\K_n \stackrel{\sim}{\to}  \{\K_n\otimes^\L\Z/p^t\}_t.$$
\end{proposition}

\proof
1) This is proven in \cite[Lemma 2]{flach}. 
\iffalse
It suffices to show that the induced maps
$${}_{p^t}i^*\K\to \{ {}_{p^t}\K_n\}_n,\qquad 
i^*\K/p^t\to \{\K_n/p^t\}_n$$
are pro-isomorphisms. By induction on $t$ and the kernel-cokernel sequence 
for the composition 
$A\stackrel{p^i}{\to}A\stackrel{p}{\to}A $ we can reduce this to the case $t=1$.
\fi

2) The pro-system $\{\K_n\otimes^\L\Z/p^t\}_t$ is quasi-isomorphic to the 
pro-complex
$\{\K_n \stackrel{p^t}{\to} \K_n\}_t$, where the transition maps
in the left system are multiplication by $p$. 
It suffices to show that the left-system is Artin-Rees zero. 
For this we fix $s$ such that $p^s\geq n$ and 
show that the $(n+s)$-fold transition map in the system is the zero map.

The stalks of $\K_n$ are sections over strictly local $\Z/p^n$-algebras, and
every section of $\K_n$ over a local $\Z/p^n$-algebra
$A$ can be written as $1+x$ with $x\in pA$. It suffices to show that 
$(1+x)^{p^{n+s}}=1$ for all 
$x\in pA$. The monomials $\binom{p^{n+s}}{j}x^j$  vanish for 
$j\geq p^s\geq n$ because
$x\in pA$ and $p^n=0$ in $A$, and they vanish for $0<j<p^s$ 
by the following lemma.
\proofend

\begin{lemma}
For fixed $s$ we have $v_p(\binom{p^{n+s}}{u})>n$ for all $0<u<p^s$.
\end{lemma}

\proof
From Legendre's formula we get
$$v_p(\binom{z}{u})=\sum_{i=1}^\infty
\big( \left\lfloor{z \atop p^i}\right\rfloor-
\left\lfloor{u\atop p^i}\right\rfloor-
\left\lfloor{z-u\atop p^i}\right\rfloor\big),$$
where for a real number $x$, $\lfloor x\rfloor$ denotes the largest integer which is not larger than $x$.
For real numbers $x$ and $y$ we have 
$\left\lfloor x+y\right\rfloor-\left\lfloor x\right\rfloor-
\left\lfloor y\right\rfloor\geq 0$,
and $\left\lfloor x+y\right\rfloor-\left\lfloor x\right\rfloor-
\left\lfloor y\right\rfloor=1$,
if $x+y$ is an integer but $x,y$ are not.
If $z=p^{n+s}$ and $0<u<p^s$, then $\frac{z}{p^i}$ is an integer but 
$\frac{u}{p^i}$ is not for $i=s,\ldots ,n+s$, and the Lemma follows.
\proofend

\iffalse
The proposition implies that we can rewrite \eqref{seqc} as follows:
%\begin{theorem}\label{mainth}
%There is an exact sequence
\begin{multline}\label{complseq}
0\to H^1_\cont(\X,\K^\wedge) \to \Pic(\X)\to\Pic(\X_s) \to \\
H^2_\cont(\X,\K^\wedge) \to  H^2_\cont(\X_\bullet, \G_m)\to 
H^2_\et(\X_s,\G_m)\to H^3_\cont(\X,\K^\wedge).
\end{multline}
%\end{theorem}
\fi

%%%%%%%%%%%%%%%%%%%%%%%%%%%%%%%%%%%%%%%%%%

\section{Comparison of sequences}
The natural maps $i^*\G_{m,\X} \to i_n^*\G_{m,\X_n}$ induce maps 
$H^i_\et(\X,\G_m)\to  H^i_\cont(\X_s, \G_{m,\bullet})$
and $H^i_\et(\X,\K)\to  H^i_\cont(\X_s, \K_\bullet)$, 
hence 
%The morphisms $\X_n \to\X$ induce 
a map between the sequences \eqref{seq} and \eqref{seqc}.  
By Proposition \ref{mainth} we obtain  a diagram
\begin{equation}\label{bigdia}
\begin{CD}
0@>>> H^1_\et(\X,\K)@>>>\Pic(\X)@>>> 
\Pic(\X_s) @>>> H^2_\et(\X,\K)@>>> \\
@.@VVV @|@| @VaVV\\
0@>>> H^1_\cont(\X,\K^\wedge)@>>> \Pic(\X)@>>> 
\Pic(\X_s) @>>>H^2_\cont(\X,\K^\wedge)@>>> 
\end{CD}$$
$$\begin{CD}
@>>> \Br(\X)@>>>H^2_\et(\X_s,\G_m)@>u>>
H^3_\et(\X,\K) \\
@. @VbVV @|@VcVV \\
@>>> 
 H^2_\cont(\X_s, \G_{m,\bullet})@>>>H^2_\et(\X_s,\G_m)@>v>>H^3_\cont(\X,\K^\wedge).
\end{CD}\end{equation}
%We note that $v$ has finite image by the Proposition (or because the source
%is a torsion group and the target a finitely generated $\Z_p$-module).
%Hence $c$ also has finite image. 
\iffalse alt
$$\begin{CD}
\Pic(\X)@>>> 
\Pic(\X_s) @>>> H^2_\et(\X,\K)@>>> \Br(\X)@>>>H^2_\et(\X_s,\G_m)@>u>>
H^3_\et(\X,\K) \\
@| @VaVV @VbVV @|@VcVV \\
%\Pic(\X)@>>> 
\Pic(\X_s) @>>>H^2_\cont(\X,\K^\wedge)@>>> 
 H^2_\cont(\X_\bullet, \G_m)@>>>H^2_\et(\X_s,\G_m)@>v>>H^3_\cont(\X,\K^\wedge).
\end{CD}$$
\fi

Thus $H^1_\et(\X,\K)\cong H^1_\cont(\X,\K^\wedge)$, we have an
isomorphism $\ker a\cong \ker b$, and a sequence
$$0\to \coker a\to \coker b\to \ker c\cap \im u\to 0.$$
%Recall that $N=\Pic(\X_s)/\im \Pic(\X)$.

\begin{lemma}
We have $T_pH^2_\et(\X,\K)=0$, and the composition 
$$N\to H^2_\et(\X,\K)\to \lim H^2_\et(\X,\K)/p^r $$ 
is injective.
\end{lemma}

\proof 
Since ${}_{p^r}H^i_\et(\X,\K)$ is finite for $i\leq 2$, 
we obtain a short exact sequence 
$$ 0\to \lim H^i_\et(\X,\K)/p^r \to H^i_\cont(\X,\K^\wedge)\to 
T_p H^{i+1}_\et(\X,\K)\to 0.$$
For $i=1$, $H^1_\et(\X,\K)\cong H^1_\cont(\X,\K^\wedge)$
implies that this group is $p$-adically complete and that 
$T_p H^2_\et(\X,\K)$ vanishes. For $i=2$, the sequence implies that 
$\lim H^2_\et(\X,\K)/p^r \subseteq H^2_\cont(\X,\K^\wedge)$, 
and that the map from $N$ to the latter group is injective by 
diagram \eqref{bigdia}.
\endproof

\begin{proposition}\label{h2k}
The torsion subgroup $\Tor H^2_\et(\X,\K)$ is a finite $p$-group, and 
$H^2_\et(\X,\K)/\tor$ is an extension of $\Q^t$ by $\Zp^s$ with $s+t=r$.
\end{proposition}

\proof
Since $C=H^2_\et(\X,\K)$ is a $\Zp$-module,  
$C_\tor$ consists of $p$-power torsion. But every torsion group
contains a basic subgroup \cite[Thm.\ 10.36]{rotman}, i.e., 
$C$ has a pure subgroup $B$ which is a
direct sum of cyclic groups and such that $C/B$ is divisible. 
Since  $T_pC$ vanishes, we have $C=B$, and then the finiteness 
of ${}_pC$ implies finiteness of $C_\tor$. 

Now $\bar C=C/C_\tor$ is a $\Zp$-submodule of 
$C_\Q\cong N_\Q\cong \Q^r$ because $\Br(\X)$ is torsion. 
It thus suffices to prove the following lemma.

\begin{lemma}\label{sundt}
Every $\Zp$-submodule $M$ of $\Q^r$ with $M_\Q\cong \Q^r$ is an extension
of $\Q^t$ by $\Zp^s$ with $s+t=r$.
\end{lemma}

\proof
We proceed by induction on $r$. If $r=1$, then $M$ is a $\Zp$-
submodule of $\Q$.
Every element of $M$ can be written as $ap^u$, where $a\in \Zp^\times$ and $u\in \Z$.
If there exist elements with arbitrary large negative $u$, then $M=\Q$,
and if not, $M= p^{-v}\Zp$ for some $v$, hence $M$ is 
isomorphic to $\Zp$. If $r>1$, let $\pi:\Q^r\to \Q$ be
a non-trivial homomorphism with kernel $\Q^{r-1}$, and consider the diagram
$$\begin{CD}
0@>>> M\cap \Q^{r-1}@>>> M @>>> \pi(M)@>>> 0\\
@. @VVV @VVV @VVV \\
0@>>> \Q^{r-1}@>>> \Q^r @>>> \Q@>>> 0.
\end{CD}$$
By induction hypothesis $M\cap \Q^{r-1}$ has a free $\Zp$-submodule
 with uniquely 
divisible quotient, and $\pi(M)$ is either isomorphic to $\Zp$ or to $\Q$.
In the former case, we have $M\cong M\cap \Q^{r-1}\oplus \Zp$,
and in the latter case, $M$ still has a free $\Zp$-submodule with uniquely 
divisible quotient.
\proofend

\begin{corollary}\label{br1}
We have 
$$ \ker\left( \Br(\X)\to \Br(\X_s)\right) \cong 
(\Q/\Z')^s\oplus(\Q/\Z)^t\oplus P$$
with $s+t=r$ and $P$ a finite $p$-group.
Moreover $s=0$ is equivalent to $r=0$.%, i.e. $\rho_{\X_s}=\rho_{\X}$.
\end{corollary}

Note that $\Br(\X_s)$ is finite by Artin's conjecture.

\proof
Since $\Br(\X)$ is torsion and $\Br(\X_s)$ is the torsion subgroup of
$H^2_\et(\X_s,\G_m)$, the kernel of
$\Br(\X)\to \Br(\X_s)$ is isomorphic to $H^2_\et(\X,\K)/N$.
Let $A$ be the kernel of the composition
$N\hookrightarrow H^2_\et(\X,\K)\to \Q^t$ of Proposition \ref{h2k}, 
and let $B$ be its image. By Proposition \ref{h2k} we obtain a diagram
$$ \begin{CD}
0@>>> A@>>> N@>>> B@>>> 0\\
@.@VVV @VVV @VVV \\
0@>>> (\Zp)^s\oplus F @>>> H^2_\et(\X,\K)@>>> \Q^t@>>>0 
\end{CD}$$
where $F$ is a finite $p$-group and all vertical maps are injective. 
Tensoring with $\Q$ we see
that $A$ is finitely generated of rank $s$ and $B$ is finitely generated
of rank $t$.  We obtain a short exact sequence of cokernels
$$ 0\to (\Q/\Z')^s\oplus F'\to H^2_\et(\X,\K)/N\to (\Q/\Z)^t\to 0,$$
where $F'$ is a quotient of $F$. Because $(\Q/\Z')^s$ is injective, 
$H^2_\et(\X,\K)/N$ is an extension of 
$(\Q/\Z)^t\oplus(\Q/\Z')^s$ by $F'$, and 
this is isomorphic to $(\Q/\Z)^t\oplus(\Q/\Z')^s\oplus \tilde F$ with
$\tilde F$ a quotient of $F'$ by the following Lemma.
\proofend

%But as an extension of torsion cotorsion groups, $E$ is a torsion cotorsion group,
%hence a direct sum of a bounded group and a divisible group (Fuchs, Cor. 54.4).
%The coranks can then be calculated by considering the $l^r$-torsion subgroups for
%all $l$ and increasing $r$. 

\begin{lemma}
Let $E$ be an extension of a divisible group $D$ by a finite group $F$.
Then $E\cong  D\oplus F'$ for $F'$ a quotient of $F$.
\end{lemma}

\proof
Let $E'$ be the maximal divisible subgroup of $E$, and let $K$ and $D'$
be the kernel and image of the composition $E'\to E\to D$, respectively,
so that we obtain a diagram with vertical inclusions
$$\begin{CD}
0@>>> K@>>>E' @>>>D' @>>> 0\\
@.@VVV @VVV @VVV\\
0@>>> F@>>>E @>>>D @>>> 0.
\end{CD}$$
Finiteness of $F$ implies finiteness of $K$, which implies that 
the divisible groups $E'$, $D'$, and $D$ all have the same $l$-corank
for all $l$.  Hence the injection $D'\to D$
is an isomorphism, and $E=E'\oplus F/K$.
\proofend

%%%%%%%%%%%%%%%%%%%%%%%%%%%%%%%%%%%%%%%%

\section{Relationship to the chern class map, examples}
%\section{Abelian and K3-surfaces}
We have the following theorem
of Flach-Siebel \cite[Lemma 1]{flach}.

\begin{theorem}\label{flsi}
We have $H^i(X,\O_X)\cong 
H^i_\cont(\X,\K^\wedge)\otimes_{\Z_p}\Q_p$. 
\end{theorem}

The Theorem together with the sequence \eqref{seqc}
and Proposition \ref{mainth}
%\begin{multline*}
%0\to H^1_\cont(\X_s,\K_\bullet) \to \Pic(\X) \to\Pic(\X_s)\\ 
%\to H^2_\cont(\X_s,\K_\bullet)\to H^2_\cont(\X_s,\G_{m,\bullet}) \to 
%H^2_\et(\X_s,\G_m).
%\end{multline*}
induces an injection 
$$\beta: N_\Q\to H^2_\cont(\X_s,\K_\bullet)_\Q\cong H^2(X,\O_X)$$
which we can extend to a map
$$\beta_{\Q_p}: N\otimes\Q_p\to H^2(X,\O_X).$$

\begin{theorem}\label{chern}
We have $s=\dim_{\Q_p}\im \beta_{\Q_p}$ and $t=\dim_{\Q_p}\ker \beta_{\Q_p}$.
\end{theorem}

In other words, 
$s$ is the dimension of the $\Q_p$-vector space spanned
by the abelian group $N$ of rank $r$.

\proof
%By Theorem \ref{flsi}, the sequence \eqref{seqc} takes, up to finite groups, 
%the form
%\begin{equation}\label{mains}
%0\to H^1(\X,\O_\X)\to \Pic(\X) \to \Pic(\X_s) \to H^2(\X,\O_\X)
%\to H^2_\cont(\X_s,\G_{m,\bullet})\to H^2_\et(\X_s,\G_m).
%\end{equation}
Consider the following commutative diagram. 
$$\begin{CD}
@.0@.0@.0@.0\\
@.@VVV @VVV@VVV@VVV\\
0@>>> H^1_\cont(\X_s,\K_\bullet)@>>> \Pic\X^{\wedge p} @>f>> 
\Pic\X_s^{\wedge p}@>g>> H^2_\cont(\X_s,\K_\bullet)\\
%@VVV @VVV @| @VVV \\
%H^1(\X,p\O)@>>> \hat H^2(\X,\Z(1)) @>>> H^2(\X_s,\Z(1)) @>>> H^2(\X,p\O)\\
@.@| @VVV @VVV @| \\
0@>>> H^1_\cont(\X_s,\K_\bullet)@>>> H^2_\et(\X,\Z_p(1)) @>>> 
H^2_\et(\X_s,\Z_p(1)) @>>> H^2_\cont(\X_s,\K_\bullet)\\
@.@VVV @VVV @VVV @VVV\\
@.0@>>> T_p\Br\X @>\alpha'>> T_p\Br\X_s@>>>0  \\
@.@.@VVV@VVV\\
@.@.0@.0
\end{CD}$$
The upper (non-exact) row 
is obtained by completing the cohomology groups in \eqref{seqc}, 
the exact middle row is obtained by $p$-completing the coefficients 
in \eqref{complcoh}.
%\begin{equation}\label{complcoh}
%0\to \K_n \to i_n^*\G_{m,\X_n} \to \G_{m,\X_s}\to 0,
%\end{equation}
The columns are exact coefficient sequences.
The middle left horizontal map is injective because 
$H^1_\et(\X,\Z_p(1)) \to  H^1_\et(\X_s,\Z_p(1))$ is surjective, 
hence so is the upper left horizontal map. 
A diagram chase shows that $\ker g/\im f\cong  \ker \alpha'$
(which has rank $t$ by Theorem \ref{maintheo}). 
On the other hand, the diagram
$$\begin{CD}
\Pic\X\otimes\Z_p @>>> \Pic\X_s\otimes\Z_p@>>>N\otimes\Z_p@>>> 0\\
@VVV @| @| \\
\Pic\X^{\wedge p} @>f>> \Pic\X_s^{\wedge p}@>>>
N^{\wedge p}@>>> 0\\
\end{CD}$$
shows that $\Pic\X_s^{\wedge p}/\im f\cong N\otimes\Z_p$, a $\Z_p$-module
of rank $r$, and that $(\im g)\otimes_{\Z_p}\Q_p=\im\beta_{\Q_p}$.
Combining this with the canonical short exact sequence
$$ 0\to \ker g/\im f\to \Pic\X_s^\wedge/\im f\to 
\Pic\X_s^\wedge/\ker g\to 0.$$
%Since completion is exact on finitely generated groups, 
%$\Pic\X_s^\wedge/\im f$ has rank $r$ over $\Z_p$. 
we see that
$\Pic\X_s^{\wedge p}/\ker g\cong \im g$ has rank $s$.
\qed
%If we consider the same diagram at $l$ instead of $p$,
%and replace all the outer groups by $0$,
%then the proper base change shows that $\ker \alpha'$ has rank $r$.

\medskip

%\begin{corollary}
%In the good reduction case, the image of the map 
%$\Pic\X_s^\wedge \to H^2(X,\O_X)$ has rank $s$.
%\end{corollary}

If $W$ denotes the Witt vectors of the residue field, then 
in the good reduction case we have a commutative diagram by the work 
of Berthelot-Ogus \cite[Cor. 3.7]{bo}, 
%and that the diagonal composition agrees with the map $g$
$$\begin{CD}
\Pic\X @>c_{dR}>> H^2(X,F^1\Omega_X^\bullet) @>>> 
H^2(X,\Omega_X^\bullet)
@>>> H^2(X,\O_X)\\
@VVV @. @| \\
\Pic\X_s @>c_{cris}>> H^2_{crys}(\X_s/W)^{F=p} @>>> 
H^2_{crys}(\X_s/W)\otimes_WK,
\end{CD}$$
and we obtain $\beta_\Q$ because the upper row is the zero-map.

\subsection*{Examples: Abelian and K3-surfaces}
\iffalse
Recall the exact sequence of Proposition \ref{h2k},
$$0\to \Zp^s\oplus (finite) \to H^2_\et(\X,\K) \to \Q^t \to 0$$ 
with $s+t=r$ and $s=0$ implies $t=0$.

\begin{corollary}\label{br2}
We have $s+\rk T_pH^3_\et(\X,\K)=f\cdot h^{0,2}$. 
%In particular, $s\leq f\cdot h^{0,2}$.
In particular, if $h^{0,2}=0$, then $\ker (\Br(\X)\to \Br(\X_s))$ is 
finite and $r=0$.
\end{corollary}

\proof 
We have an exact sequence
$$ 0\to H^2_\et(\X,\K)^{\wedge p} \to H^2_\cont(\X,\K^\wedge)\to 
T_p H^3_\et(\X,\K)\to 0,$$
and Proposition \ref{h2k}  implies that the left term is $\Z_p^s$.
The final statement follows from the Theorem by using Corollary \ref{br1}.
\endproof
\fi

We calculate $\Br(\X)$ for $\X$ an abelian scheme or  
a family of K3-surfaces over $\Z_p$.

\begin{theorem}\label{k3ab}
Let $\X$ be an abelian scheme or a family of K3-surfaces over $\Z_p$.
If $r=0$, then $\Br(\X)$ is finite. If $r>0$, then 
$$\Br(\X)\cong (\Q/\Z')\oplus (\Q/\Z)^{r-1}\oplus (finite).$$
\end{theorem}

\proof
Since Tate's conjecture is known for abelian varieties \cite{tate}
and K3-surfaces over a finite field \cite{madapusi}, \cite{k3char2}, 
we know that $\Br(\X_s)$ is finite. Then $r=0$ implies that $\Br(\X)$ is
finite.
If $r>0$, then since $H^2(\X,\O_\X)\cong\O$ we obtain that $s\leq 1$
by Theorem \ref{chern}, but 
since $s=0$ implies $r=s+t=0$ we must have $s=1$.
\proofend

We give an example of an abelian surfaces with $r>1$, showing that 
$\Br(\X)$ contains a divisible $p$-group.

\begin{proposition}
If $\X_s$ is a simple abelian surface over $\F_p$, then $\rk \Pic(\X_s)=2$.
If $\X_s$ is the product of two elliptic curves $E_1$ and $E_2$ over $\F_p$,
then $\rk \Pic(\X_s)=4$ if $E_1$ are $E_2$ are isogenous, and 
$\rk \Pic(\X_s)=2$ if they are not.
\end{proposition}

\proof
This follows by considering Weil-numbers.
\proofend

The Picard numbers of $\XK$ can be calculated explicitly in many cases.

\begin{example}
Let $\X/\Z_p$ be the Jacobian of the smooth projective curve of genus two
defined by the equation $y^2 = x^5-1$. 
The $5$th roots of unity $\mu_5$ act on $\X$ in the obvious way, so that 
$\Z[\zeta_5]$ acts on any factor of $\XK_{\bar K}$. 
By the classification of endomorphism
algebras of abelian varieties we see that $\XK_{\bar K}$ is simple
and that  $\End(\XK_{\bar K} )_\Q=\Q(\zeta_5)$. 
The rank of the N\'eron-Severi group of $\XK$ is $1$,
because it is a subgroup of $\End(\XK )_\Q$ which is $\Q$ 
by \cite{page}. 

If $p\equiv -1 \pmod 5$, 
then $\X_s$ has good reduction at $p$ and $\X_s$ is isogenous to 
$E^2$, where $E$ is an elliptic curve over $\F_p$ satisfying
$|E(\F_p)| = p+1$. Hence the rank of the N\'eron-Severi group of the special
fiber is $4$,  % and the rank of the N\'eron-Severi group of $\X$ is $2$,
$r=3$, and we obtain 
$$ \Br(\X)\cong (\Q/\Z')\oplus (\Q/\Z)^2\oplus (finite).$$
\end{example}

\section{The inverse system of Brauer groups}
We discuss the maps in the diagram 
\begin{equation}\label{H2} 
\begin{CD}
@.@. \Br(\X)@>>>H^2_\et(\X_s, \G_m)\\
@.@.@VbVV@AAA\\
0@>>> \lim^1 \Pic(\X_n)@>>> H^2_\cont(\X_s, \G_{m,\bullet})
@>>> \lim H^2_\et(\X_n,\G_m)@>>> 0,
\end{CD}
\end{equation}
where the lower sequence is \eqref{slim} for $i=2$. As a first
result we have

\begin{proposition}\label{limbrauer}
The kernel of $\lim H^2_\et(\X_n,\G_m)\to H^2_\et(\X_s,\G_m)$ is a 
finitely generated $\Z_p$-module of rank at most $f\cdot h^{0,2}$.
\end{proposition}

\proof
From the exact sequences
\begin{equation}\label{ascent}
H^2(\X_s,\O_{\X_s})\to H^2_\et(\X_n,\G_m)\to H^2_\et(\X_{n-1},\G_m)\to H^3(\X_s,\O_{\X_s})
\end{equation}
we inductively see that the kernel $K_n$ %and cokernel $C_n$ 
of $H^2_\et(\X_n,\G_m)\to H^2_\et(\X_s,\G_m)$ is a finite $p$-group. 
Taking the limit we obtain 
$$ 0\to \lim K_n \to \lim H^2_\et(\X_n,\G_m)\to H^2_\et(\X_s,\G_m).$$
This shows that the kernel is a pro-$p$ group.
But the finitely generated $\Z_p$-module $H^2(\X,\O_{\X})$ surjects onto 
the kernel of the composition
$$H^2_\cont(\X_s,\G_{m,\bullet})\to\lim H^2_\et(\X_n,\G_m)\to 
H^2_\et(\X_s,\G_m),$$ 
hence it surjects onto the kernel of the second map because the first map is
surjective.
\proofend

%We are going to analyze the situation more carefully.
%We first consider the map $a$. 
%Since $N$ has finite image in $C_\tor$, we know that 
%$N \to \bar C\to \lim \bar C/p^r$ has finite kernel. 

\begin{theorem}
The map $\Br(\X)\stackrel{b}{\longrightarrow} 
H^2_\cont(\X_s, \G_{m,\bullet})$ is injective.
\end{theorem}

Grothendieck \cite[Lemma 3.3]{brauerIII} showed that the natural map 
$\Br(\X)\to \lim \Br(\X_n)$ is injective if the system $(\Pic(\X_n))_n$ is 
Mittag-Leffler. But
if $(\Pic(\X_n))_n$ is Mittag-Leffler, then 
$H^2_\cont(\X_s, \G_{m,\bullet})\cong \lim H^2_\et(\X_n,\G_m)$, hence 
the Theorem is a generalization \cite[Lemma 3.3]{brauerIII}. The theorem
follows by diagram \eqref{bigdia} from the following proposition.

\begin{proposition}
The map $a:H^2_\et(\X,\K)\to H^2_\cont(\X,\K^\wedge) $ is injective, and 
its cokernel is the finitely generated free $\Z_p$-module 
$T_p H^3_\et(\X,\K)$ if $r=0$, and the direct sum of 
$T_p H^3_\et(\X,\K)$ and an uncountable, uniquely divisible group
if $r>0$. In particular, the cokernel of $a$ is a torsion 
free $\Zp$-module.
\end{proposition}

\proof
%Let $C =H^2_\et(\X,\K) $. 
We showed that $C =H^2_\et(\X,\K) $ is a $\Zp$-module with finite torsion
and $\bar C=C/C_\tor$ is an extension of $\Q^t$ by $\Zp^s$ with $s+t=r$.
Since $\Br(\X)$ is torsion, 
$$ \ker a%\big( H^2_\et(\X,\K)\to H^2_\cont(\X,\K^\wedge)\big)
\cong 
\ker\big(\Br(\X)\to H^2_\cont(\X_\bullet, \G_m)\big)\subseteq C_\tor$$
%$\ker a\cong \ker b\subseteq C_\tor$
is finite. Moreover, $a$ factors as 
$$a :C\to C^\wedge \to H^2_\cont(\X,\K^\wedge),$$
where the second map is injective with cokernel $T_p H^3_\et(\X,\K)$.
We have a diagram
$$ \begin{CD}
0@>>> C_\tor @>>> C@>>> \bar C@>>> 0\\
@.@| @VVV@VVV \\
0@>>> C_\tor^\wedge @>>> C^\wedge @>>> \bar C^\wedge @>>> 0.
\end{CD}$$
The lower row is exact on the left because $\bar C$ is torsion free and
on the right because $\lim^1C_\tor/p^r=0$. The injectivity of 
$C_\tor \stackrel{\sim}{\to} C_\tor^\wedge \to C^\wedge$ implies
%$\ker a\subseteq C_\tor$ implies 
that $a$ is injective.
By Corollary \ref{compllemma}, the cokernel is the direct sum
of $T_p H^3_\et(\X,\K)$ and the uniquely divisible group 
$\lim^1 (H^2_\et(\X,\K),p)$. By Proposition \ref{h2k}, and the fact that 
$\lim(\Zp,p)=\lim^1 (\Q^t,p)=0$, we have a sequence
$$ 0\to \lim (H^2_\et(\X,\K),p)\to \Q^t\to \lim^1(\Zp^s,p)
\to \lim^1 (H^2_\et(\X,\K),p)\to 0. $$
Taking the long exact derived lim-sequence of the sequence of inverse systems 
$$ 0\to (\Zp^s,p)\stackrel{(p^n)}{\to} (\Zp^s,\id)\to (\Z/p^n)^s\to 0$$
we obtain $(\Z_p/\Zp)^s\cong \lim^1(\Zp^s,p)$, hence the result.
\proofend

\begin{corollary}
Assuming finiteness of $\Br(\X_s)$, $\Br(\X)$ agrees with the torsion
subgroup of $\lim^1\Pic(\X_n)$ up to finite groups. 
\end{corollary}

\proof
Finiteness of $\Br(\X_s)$ implies finiteness of  $\Br(\X_n)$ 
by the sequence \eqref{ascent}, and this implies that 
$\lim H^2_\et(\X_n,\G_m)$ 
is a pro-finite group. It follows that any divisible group maps to zero in
$\lim H^2_\et(\X_n,\G_m)$, hence the divisible part of $\Br(\X)$ 
injects into $\lim^1\Pic(X_n)$. On the other hand, 
the torsion free group $\coker a$ is a subgroup of the cokernel of 
$\Br(\X)\to H^2_\cont(\X_\bullet, \G_m)$ of finite index by 
\eqref{bigdia} and finiteness of $\Br(\X_s)$.
\proofend

We are comparing our results to results about $\lim^1 \Pic(\X_n)$. 
The groups $\Pic(\X_n)$ are finitely generated of 
constant rank, and
as the derived limit of finite groups vanishes, we can 
consider the torsion free quotients $\overline{ \Pic(\X_n)}$ instead. 
Since $\Pic(\X_n)\to \Pic(\X_{n-1})$ has finite kernel and cokernel,
the maps $\overline{ \Pic(\X_n)}\to \overline{ \Pic(\X_{n-1})}$ are injective,
hence the images $T$ of $\Pic(\X)$ in each group are isomorphic. 
We obtain an exact sequence of pro-systems
$$
0\to T \to \overline{ \Pic(\X_n)} \to Q_n\to 0,
$$
hence $\lim^1 \Pic(\X_n)\cong\lim^1\overline{ \Pic(\X_n)}\cong \lim^1 Q_n$, 
where each $Q_n$ is finitely generated of rank $r$.

By Jensen \cite[Thms. 2.5, 2.7]{jensen}, 
if the groups $A_i$ in a countable 
pro-system are finitely generated, then $\lim^1A_i\cong \Ext(M,\Z)$, where 
$M=\colim \Hom(A_i,\Z)$ is a countable torsion free group. 
Moreover, if $\lim^1A_i$ does not vanish, then 
\begin{equation}\label{jenn}
\lim^1A_i\cong  \Q^{n_0}\oplus \bigoplus_p (\Q_p/\Z_p)^{n_p} 
\end{equation}
where $n_0$ is the cardinality of the continuum $2^{\aleph_0}$, and 
$n_p$ is either $2^{\aleph_0}$ or finite (possibly zero). 
We give a more precise statement in a special situation.
%for $A_i$ with constant rank 
%$r$ and finite cokernels of transition maps. 

\begin{theorem}\label{thm8.1}
Let $A_i$ be an inverse system of finitely generated groups of constant 
rank $r$ and transition maps with finite cokernel. 
If $\lim^1A_i$ does not vanish, then 
%$$\lim^1A_i\cong  \Q^{n_0}\oplus \bigoplus_p (\Q_p/\Z_p)^{n_p} $$
%with $n_0$ the cardinality of the continuum and 
$0\leq n_p\leq r$ in \eqref{jenn}.
If the cokernels of the maps in the system are finite $p$-groups, 
then $\lim^1A_i$ vanishes or 
$n_p<n_l$ for all $l\not=p$, and $n_l=n_{l'}$ for $l,l'\not=p$.
\end{theorem}

\proof
We can assume that each group $A_i$ is a free abelian group of rank $r$
and proceed by induction on $r$. If $r=1$, we 
let $M=\colim \Hom(A_i,\Z)$. Choosing
any non-zero element of $M$ identifies $M\otimes\Q$ with $\Q$, hence
the inclusion $M\to M\otimes\Q\cong\Q$ identifies $M$ with a
subgroup of $\Q$, which is of the form $\Z[\{p^{-e_p}\}_p]$, 
where $p$ runs through 
the primes and $e_p$ is an integer or infinity. 
A different choice of an element of $M$ changes finitely many $e_p$
by a finite amount.

\begin{lemma}\label{rank1}
Let $M\cong \Z[\{p^{-e_p}\}_p]\subseteq \Q$, where $p$ runs through 
the primes and $e_p$ is a non-negative integer or infinity. 
If  all $e_i$ are finite and almost all vanish, then 
$M\cong \Z$ and $\Ext(M,\Z)=0$. Otherwise 
$$\Ext(M,\Z)\cong \Q^{n_0}\oplus\bigoplus_p (\Q_p/\Z_p)^{n_p},$$
where $n_0=2^{\aleph_0}$, $n_p=0$ if $e_p$ is infinity, and $n_p=1$ if $e_p$ is finite.
\end{lemma}

For example, 
$$
\Tor \Ext(M,\Z)=\begin{cases}
(\Q/\Z)', &M=\Z[p^{-\infty}] \\
\Q_p/\Z_p,  &M=\Z_{(p)}\\
\Q/\Z,&M=\Z[p^{-1}|\text{infinitely many } p]
\end{cases}$$

\proof
The case $M\cong\Z$ is easy, so that we assume $M\not\cong \Z$.
The long exact $\Ext^i(-,\Z)$ sequence associated to the short
exact sequence $0\to \Z\to M\to M/\Z\to 0$ together with $\Hom(M,\Z)=0$
because $M\not\cong \Z$ gives
$$0\to \Z\to \Ext(M/\Z,\Z)\to  \Ext(M,\Z)\to 0.$$
Let us first consider the torsion subgroup $\Tor \Ext(M,\Z)$. 
The six term sequence associated to derived tensor product 
$-\otimes\Z/p^r$ together
with the fact that $\Ext(M,\Z)$ is divisible gives
\begin{equation}\label{extm}
0\to  {}_{p^r} \Ext(M/\Z,\Z)\to {}_{p^r} \Ext(M,\Z)\to \Z/p^r\Z\to 
\Ext(M/\Z,\Z)/p^r\to 0.
\end{equation}
Now $M/\Z \cong \bigoplus_p \Z/p^{e_p}$, where we set
$\Z/p^{e_p}=\Q_p/\Z_p$ if $e_p$ is infinity, and then 
\begin{equation}\label{yuio}
\Ext(M/\Z,\Z)\cong\Hom(M/\Z,\Q/\Z)\cong \prod_p \Z_p/p^{e_p},
\end{equation}
where we set $p^{e_p}=0$ if $e_p$ is infinity.
In particular, if $e_p$ is finite, then the left and right groups
in \eqref{extm} have the same cardinality,  
so that ${}_{p^r} \Ext(M,\Z)$ has cardinality $p^r$ for all $r$, 
hence the $p$-primary torsion of $\Ext(M,\Z)$ is $\Q_p/\Z_p$. 
On the other hand, 
if $e_p$ is infinity, then the left group in \eqref{extm} vanishes
whereas the right group is isomorphic to $\Z_p/p^r\Z_p$, 
hence we obtain ${}_{p^r} \Ext(M,\Z)=0$.

Finally, note that $\Ext(M/\Z,\Z)$, hence $\Ext(M,\Z)$,  
has the same cardinality as the continuum by our hypothesis on $e_p$
and the description in \eqref{yuio}.
Since $\Ext(M,\Z)$ is divisible and $\Tor \Ext(M,\Z)$ is countable,
we get the statement of the Proposition.
\proofend

%\proof(Theorem \ref{thm8.1})
We continue the proof of Theorem \ref{thm8.1}.
The case $r=1$ follows from the Lemma 
because $\lim^1A_i\cong \Ext(M,\Z)$ for $M=\colim \Hom(A_i,\Z)$. 
If the transition maps have $p$-groups as the cokernel, then $M\cong \Z$ 
or $M\cong \Z[p^{-\infty}]$ in which
case we get the claimed statement on the $n_p$.

For general $r$, we can, by performing elementary column operations 
(which corresponds to changing the basis of the next group in the 
inverse system), assume that the transition maps are given by matrices 
$M_i=
\begin{pmatrix}a_i&0&0\\ \star &\star&\star\\ 
\star &\star&\star\end{pmatrix}$.
Thus there is a subsystem $(A'_i)$ consisting of free groups of rank $r-1$ 
and a quotient system $(A''_i)$ consisting of free groups of rank $1$.
By hypothesis $a_i\not=0$. 
If the $a_i$ are $\pm1$ for almost all $i$, then  
$\lim A''_i\cong \Z$ and $\lim^1 A''_i=0$,
and we have a sequence
$$ 0\to \lim A'_i \to  \lim A_i \to  \Z \stackrel{\delta}{\to}  \lim^1 A'_i \to  \lim^1 A_i \to 0.$$ 
If $\delta$ has finite image, then the parameters $n_p$ for $(A_i)$ and $(A'_i)$
agree. If $\delta$ has infinite
image, then the parameters $n_p$ of $A_i$ are one larger than the 
parameters for $A'_i$.
If $a_i$ is different from $\pm1$ for infinitely many $i$, then  
$\lim A''_i=0$ and we obtain a sequence
$$ 0\to  \lim^1 A'_i \to  \lim^1 A_i \to  \lim^1 A''_i \to 0.$$ 
In this case the parameters $n_p$ of $A_i$ are the sum of the parameters
$n_p$ of $A'_i$ and of $A''_i$.
%The proof of the second statement is the same, using that in the base step
%$\Tor \Ext(\Z[p^{-\infty}],\Z)=(\Q/\Z)'$.
\proofend

\end{document}